\documentclass[11pt]{article}
\usepackage{amsmath, amssymb, amsfonts, verbatim, amsthm, xcolor}
\usepackage{graphicx}
\newtheorem{theorem}{Theorem}[section]

\newtheorem{proposition}[theorem]{Proposition}

\newtheorem{lemma}[theorem]{Lemma}

\def\qed{\hfill $\Box$\medskip}

\def\IC{{\mathbb C}}

\def\cl{{\bf cl}}
\def\bg{{\bf g}}
\def\la{{\langle}}
\def\ra{{\rangle}}
\def\bd{{\bf d}}
\def\bc{{\bf c}}
\def\({\left (}
\def\){\right )}

\def\diag{{\rm diag}\,}
\def\tr{{\rm tr}\,}

\def\span{{\rm span}\,}

\def\bd{{\bf d}}

\def\bx{{\bf x}}
\def\by{{\bf y}}

\def\conv{{\bf conv}}

\begin{document}
\openup 1\jot
\title{The $C$-numerical range and Unitary dilations}
\author{Chi-Kwong Li\thanks{This research was supported by the 
Simons Foundation Grant 851334.} \\
Department of Mathematics, College of William \& Mary, \\
Williamsburg, VA 23185, USA. ckli@math.wme.du}
\date{}
\maketitle

\centerline{\bf In memory of Mrs.\ Tso-Mei Au-Yeung.}

\begin{abstract}
For an $n\times n$ complex matrix $C$, 
the $C$-numerical range of a bounded linear operator
$T$ acting on a Hilbert space of dimension at least $n$ is the set
of complex numbers $\tr(CX^*TX)$, where $X$ is a partial isometry
satisfying $X^*X = I_n$. It is shown that 
$$\cl(W_C(T)) = \cap \{\cl(W_C(U)): U \hbox{ is a unitary dilation of } T\}$$
for any contraction $T$ if and only if $C$ is a rank one normal matrix.
\end{abstract}



\section{Introduction}

Let $B(H)$ be the set of bounded linear operators acting on the Hilbert space
$H$ equipped with the inner product $\la x, y\ra$. If $H$ has dimension $n$, 
then $H$ is identified with $\IC^n$ 
with the inner product $\la x, y\ra = y^*x$, and 
$B(H)$ is identified with the set $M_n$ of $n\times n$ complex matrices.
The {\it numerical range} of $T \in B(H)$ is defined as
$$W(T) = \{\langle Tx,x\rangle: x \in H, \langle x, x\rangle = 1\},$$
which is a useful concept in studying operators and matrices; 
see, e.g., \cite{Halmos1,HJ,WG}. In particular, there are interesting connections
between the study of numerical range and dilation theory; e.g., see 
\cite{CL0,CL,Halmos2}.
Recall that an operator  $\hat T \in B(K)$ 
is a {\it dilation} of $T\in B(H)$ if $K$ is a Hilbert space containing 
$H$, and 
$\hat T$ has an operator matrix of the form $\begin{pmatrix}
T & \star \cr \star & \star\cr\end{pmatrix}$ with respect to 
the space decomposition $H \oplus H^\perp$.
Equivalently, there is a partial isometry $X: H \rightarrow K$ with 
$X^*X = I_H$ and $X^*\hat TX = T$.  We also say that $T$ is a {\it compression} 
of $\hat T$.

\medskip
It is easy to show that $W(T) \subseteq W(\hat T)$ if $\hat T$
is a dilation of $T$.
Suppose  $T\in B(H)$ is a contraction, i.e., $\|T\| \le 1$. Then $T$
has a unitary dilation 
$$U = \begin{pmatrix} T & \sqrt{I - TT^*}\cr \sqrt{I-T^*T} & -T^*\cr
\end{pmatrix} \in B(H \oplus H).$$
It was conjectured in \cite{Halmos2} that for any contraction $T \in B(H)$
\begin{equation}\label{classical}
W(T) = \cap \{ W(U): U \hbox{ is a unitary dilation of } T\}.
\end{equation}
However, counter-examples for the conjecture were given in \cite{D}.
In particular, there is a normal operator $T$ with $\|T\|\le 1$ such that
equality (\ref{classical}) fails. Denote by $\cl(X)$ the closure of a set $X$ in $\IC$.
It was shown in \cite{CL} that for any contraction $T \in B(H)$,
$$\cl(W(T)) = \cap \{ \cl(W(U)): U \hbox{ is a unitary dilation of } T\}.$$
This result has been refined, and extended to other types of generalized numerical
ranges; e.g., see \cite{BGT,BT,DM,GLW}.
In particular, the authors in \cite{DM} considered the extension of 
the dilation result to the $C$-numerical range
defined as follows. 
Let $C \in M_n$ and $T\in B(H)$ with $\dim H \ge n$.
Define the $C$-{\it numerical range} of  $T$ as
$$W_C(T) = \{ \tr(CY): Y \in M_n \hbox{ is a compression of } T\}.$$
For $H$ with $\dim H \ge n$,
if we regard $C \oplus {\bf 0}$ as a finite rank operator in $B(H)$, then 
$$W_C(T) = \{\tr[(C \oplus {\bf 0})V^*TV]: V\in B(H) \hbox{ unitary}\}.$$
If $C$ is a rank one normal matrix with a nonzero eigenvalue $\gamma$, then 
$W_C(T) = \gamma W(T)$;
if  $C = I_n\in M_n$ then $W_C(T)$ reduces to the 
$n$-numerical range of $T$ consisting of complex numbers of the form
$\sum_{j=1}^n \la Tx_j, x_j\ra$ for an orthonormal set 
$\{x_1, \dots, x_n\}\subseteq  H$.
If $C$ is normal with eigenvalues $c_1,\dots, c_n$, then 
$W_C(T)$ consists of numbers of the form
$\sum_{j=1}^n c_j \la Tx_j, x_j\ra$ for an orthonormal set 
$\{x_1, \dots, x_n\}\subseteq  H$, and the set is also referred to as 
the $c$-numerical range of $T$, denoted by
$W_c(T)$, with $c = (c_1, \dots, c_n)$.
One may see \cite{GS,Li} and their references for some 
basic background of the $C$-numerical range, and how it can be
used to study matrices and operators.

In \cite{DM}, the authors considered the extension of (\ref{classical}) 
to the $C$-numerical range. Examples of $C \in M_n$ and contraction $T \in B(H)$
with $\dim H \ge n$ are given such that the following  equality fails
\label{WCA-i}
$$\cl(W_C(T)) = \cap \{ \cl(W_C(U)): U \hbox{ is a unitary dilation of } T\}.$$
In this paper, we characterize  $C \in M_n$ such that the above equality holds 
for any contraction $T \in B(H)$ with $\dim H \ge n$ by proving the following.

\begin{theorem} \label{key} 
Let $n \ge 2$ and $C \in M_n$ be nonzero. The following conditions
are equivalent. 
\begin{itemize}
\item[{\rm (a)}] The matrix $C$ is a rank one normal matrix.
\item[{\rm (b)}] For any contraction $T \in B(H)$ with $\dim H \ge n$, 
\begin{equation}\label{intersection}
\cl(W_C(T)) = \cap  \{ \cl (W_C(U)): U \hbox{ is a unitary dilation of } T \}.
\end{equation}
\item[{\rm (c)}] For any rank one nilpotent contraction $T \in M_n$, 
$$W_C(T) = \cap \{ W_C(U): U \hbox{ is a unitary dilation of }
T\}.$$
\end{itemize}
\end{theorem}

Several remarks concerning Theorem \ref{key} are in  order.

\begin{itemize}
\item Suppose $C = [\gamma] \in M_1$. Then $W_C(T) = \gamma W(T)$.
Then (\ref{intersection}) holds for all contractions
$T \in B(H)$ with $\dim H \ge 1$ by the result in \cite{CL}. So, we exclude this
case in the theorem.

\item
By Theorem \ref{key}, we see that (\ref{intersection}) holds for all contractions
$T \in B(H)$ with $\dim H \ge n$ if and only if 
$C$ is a rank one normal matrix so that $W_C(A) = \gamma W(A)$ with $\gamma = \tr C$.
Thus, the dilation result for the $C$-numerical range only admits a trivial extension.

\item
In general, $W_C(T)$ may not be convex, may not be closed, 
and may have a complicated geometrical shape. 
By condition (c),   to characterize $C\in M_n$ satisfying  
(\ref{intersection}) for all contractions in $B(H)$,
one only needs to check that
(\ref{intersection}) holds for rank one nilpotent  contractions $T$ in $M_n$,
where $W_C(T)$ is always a closed circular disk centered at the origin; see 
Lemma \ref{cor}.

\item In fact, the implication (c) $\Rightarrow$ (a)
can be strengthened to the following.
 \end{itemize}

\begin{theorem} \label{main} Suppose $n \ge 2$. If $C \in M_n$ is nonzero and is
not a rank one normal matrix, then there is a rank one nilpotent contraction 
$T \in M_n$ such that $W_C(T)$ is a closed circular disk centered at the origin
with radius $r$, and there is a positive number $d>0$ such that
$r+d \in \cl(W_C(U))$  for any unitary dilation $U$ of $T$.
Consequently, $W_C(T) = \cl(W_C(T))$ is a proper subset of 
\begin{equation}
\label{eq-fail}
\cap \{ \cl(W_C(U)): U\hbox{ a unitary dilation of } T\}.
\end{equation}
\end{theorem}
 
The proofs of Theorem \ref{key} and Theorem \ref{main}
will be presented in  the next section.  
We note that the $C$-numerical range can be defined using 
a trace class operator $C\in B(H_0)$, where $H_0$ is a separable
Hilbert space.  Our results and proofs
can be readily extended to a trace class operator
$C$ acting on a separable Hilbert space $H_0$ if
$M_n$ is replaced by $B(H_0)$  in the statements of 
Theorem \ref{key} and Theorem \ref{main}.

\section{Proofs}

We always assume that  $n \ge 2$, and denote by
$\{E_{11}, E_{12}, \dots, E_{nn}\}$  the standard basis for $M_n$.
The following properties for the $C$-numerical range are known; see \cite{Li}.

\begin{itemize}
\item  Let $C, D \in M_n$ and $S,T \in B(H)$ with $\dim H \ge n$. 
If $C = U^*DU$ and $S = V^*TV$
for some unitary $U, V$, then
$W_C(T) = W_D(S)$.

\item Let $C, A \in M_n$. Then $W_C(A) = W_A(C)$.

\item For $C \in M_n$ and $T \in B(H)$ with $\dim H \ge n$,
$W_C(\xi_1 I + \xi_2T) = \xi_1(\tr C) + \xi_2 W_C(T)$.

\item Let $C \in M_n$. If $T \in B(H)$ is a compression of $\hat T$
and $\dim H \ge n$, then $W_C(T) \subseteq W_C(\hat T)$.
\end{itemize}
We will often use the fact that $\mu \in W(T)$ for $T \in B(H)$ if and only if 
there is a unitary operator $U\in B(H)$ such that $\mu$ is the $(1,1)$ entry 
of $U^*TU$. In case $\dim H = n$
and $\mu$ is a boundary point of $W(T)$, we have the following; see \cite{HJ,Li}.

\begin{lemma} \label{lemma}
Let $A=(a_{ij}) \in M_n$. If $a_{11}$ lies on the boundary of 
$W(A)$, then  there is $\phi\in [0, 2\pi)$ such that
$e^{i\phi} a_{11} + e^{-i\phi} \bar a_{11}$ is 
the largest eigenvalue of the Hermitian matrix
$A_\phi = e^{i\phi}A + e^{-i\phi}A^*$. Consequently, 
the $(1,j)$ entry of $A_\phi$, i.e., 
$e^{i\phi} a_{1j} +  e^{-i\phi} \bar a_{j1} = 0$,
and hence
$|a_{1j}| = |a_{j1}|$ for $j = 2, \dots, n$. 
\end{lemma}

The following is known as the elliptical range theorem for the 
numerical range, see \cite{HJ,Li}. We list two special cases (a) and (b), 
which will be used frequently in our discussion.

\begin{lemma} \label{ellipse} Let $A = (a_{ij}) \in M_2$ with
eigenvalues $\lambda_1, \lambda_2$. Then $W(A)$ is an elliptical disk with
foci $\lambda_1,\lambda_2$, and length of minor axis 
$\sqrt{\tr (A^*A)-|\lambda_1|^2 -|\lambda_2|^2}$.
\begin{itemize} 
\item[{\rm (a)}] If $a_{11} = a_{22} = 0$,
then $W(A)$ is an elliptical disk with foci 
$\pm \sqrt{a_{12}a_{21}}$, and length of minor axis equal to 
$| |a_{12}|-|a_{21}| |$.
\item[{\rm (b)}] If $a_{21} = 0$, then 
$W(A)$ is an elliptical disk with foci $a_{11}, a_{22}$, and  
 length of minor axis equal to $|a_{12}|$. 
\end{itemize}
\end{lemma}

The following two results are important in our analysis;
one may see  \cite{S} and also \cite{Li,T} for their proofs. 

\begin{lemma} \label{prop1}
Let $A \in M_n$ be nonzero.
Then 
$$R = 
\min \{ \|A-\mu I\|: \mu \in \IC\}
\quad \hbox{ is equal to }
\quad \max \{|x^*Ay|: x, y \in \IC^n, \|x\| = \|y\| = 1, x^*y = 0\}.$$
Moreover,  $\mu \in \IC$ satisfies $\|A-\mu I\| = R$ if and only if
$A-\mu I$ is unitarily similar to a matrix with the $(2,1)$ entry equal to
$R = \|A-\mu I\|$, which is the only nonzero entry 
in the first column and the second row of the matrix $A-\mu I$. 
\end{lemma}

\begin{lemma}\label{cor} Let $C \in M_n$ and $T = E_{12} \in M_n$.
Then 
\begin{equation}
\label{cor:eq}
W_C(T) = \{ u^*Cv: u,v \in \IC^{n}, \|u\| = \|v\| = 1, u^*v  = 0\}
\end{equation}
is a circular disk centered at the origin with radius 
$r= \min\{ \|C-\mu I_n\|: \mu \in \IC\}.$
\end{lemma}

We have the following observation showing that for any unitary dilation 
$U$ of the rank one matrix $T= \cos\theta E_{12} \in M_n$ with $\theta \in (0, \pi/2)$,
there is  a partial isometry $X: \IC^{n+2} \rightarrow H$ such that 
$X^*X = I_{n+2}$ and $X^*UX = \hat T \in M_{n+2}$ having some specific entries.
We can then use $W_C(\hat T)$ to establish our result.

\begin{lemma} \label{dilate}
Let $T = \cos \theta E_{12} \in M_n$ and $\theta \in (0,\pi/2)$.
Suppose $U = \begin{pmatrix} T & * \cr * & * \cr \end{pmatrix} \in B(H)$ 
is a unitary dilation of  $T$. Then $\dim H \ge n+2$, and there is 
a partial isometry $X: \IC^{n+2} \rightarrow H$ such that 
$X^*X = I_{n+2}$ and 
\begin{equation}
\label{hatT}
\hat T  = X^*UX = \begin{pmatrix} T & T_{12}\cr T_{21} & T_{22}\cr\end{pmatrix}
\in M_{n+2},
\end{equation}
where 
\begin{itemize}
\item the top two rows of $\hat T$ are 
$(0, \cos \theta, 0, \dots, 0, -\sin \theta)$ and 
$(0, \dots, 0, 1,0)$;

\item rows 3 to $n$ of $\hat T$  are zero;

\item the last two rows of  $\hat T$ equal

\medskip
\qquad $(x_1, x_2 \sin \theta, x_3, \dots, x_n, 0, x_2 \cos \theta)$ \ \ and \ \ 
$(y_1, y_2 \sin\theta, y_3, \dots, y_n, 0, y_2 \cos \theta)$, 

\medskip\noindent
where
$x_4 = \cdots = x_n = 0$ if $n \ge 4$.

\end{itemize}
In particular, the 
$4\times 4$ submatrix
of $\hat T$
with row and column indices $1,2,n+1, n+2$ equals
\begin{equation}\label{matrixB}
B = \begin{pmatrix} 
0 & \cos \theta & 0 & -\sin \theta \cr
0 & 0 & 1 & 0\cr
x_1 & x_2\sin\theta & 0 & x_2\cos \theta\cr
y_1 & y_2\sin\theta & 0 & y_2\cos \theta\cr
\end{pmatrix}, \qquad \hbox{ where } 
|x_1|^2 + |x_2|^2 \le 1, \ |y_1|^2 + |y_2|^2 \le 1.
\end{equation}
\end{lemma}

\it Proof. \rm 
Since $U = \begin{pmatrix} T & \star \cr \star & \star\cr\end{pmatrix}$
is unitary, it has orthonormal rows. Thus, the first two rows of $U$ has the form
$(\overbrace{0, \cos \theta, 0, \dots, 0}^{n}, \bx)$ 
and $(\overbrace{0,\dots, 0}^n, \by)$, where
$\bx$ and $\by$ are orthogonal vectors with $\|\bx\| = \sin \theta$
and $\|\by\| = 1$. Thus, $\span\{\bx, \by\} \ge 2$ so that $\dim H \ge n+2$.
Moreover,  
there is a unitary operator 
$V = I_n \oplus V_1\in B(H)$  such that 
the first two rows of $V^*UV$ have the form 
$$(\underbrace{0, \cos\theta, 0, \dots, 0, 0, -\sin\theta}_{n+2}, {\bf 0})
\quad \hbox{ and } 
\quad (\underbrace{0, \dots, 0, 1, 0}_{n+2}, {\bf 0}).$$
Let $\hat T = \begin{pmatrix}T & T_{12} \cr T_{21} & T_{22}\cr\end{pmatrix}
\in M_{n+2}$ be the leading submatrix of  $V^*UV$.
Then the first two rows of $\hat T$
have the asserted form. Now the first two rows of $\hat T$ have unit length.
Thus, $\hat T \hat T^* = I_2 \oplus Y$ for some $Y \in M_n$.
So, all other rows of $\hat T$ are orthogonal to rows 1 and 2.
Since row 3 to row $n$ of $T$ are zero, we see that row 3 to row $n$
of $T_{12}$ must also be zero. 
Also, the last two rows of $\hat T$ must have the form
$(x_1, x_2 \sin \theta, x_3, \dots, x_n, 0, x_2\cos\theta)$ and
$(y_1, y_2 \sin \theta, y_3, \dots, y_n, 0, y_2\cos\theta)$.
If $n \ge 3$, let $W_1\in M_{n-2}$ be unitary such that 
$(x_3, \dots, x_n)W_1 = (\hat x_3, 0, \dots, 0)$. 
Let $W = I_2 \oplus W_1 \oplus I_2$. We may replace
$\hat T$ by $W^*\hat T W$ and assume that 
$x_4 = \cdots = x_n = 0$.
The assertion about $B$ is clear. \qed

\begin{lemma} \label{WCA}
Suppose $\theta \in [0, \pi/2)$, $f \in [0,1], g, x_2 \in \IC$ with $|x_2| \le 1$,
\begin{equation} \label{matrixC}
\hat C  = \begin{pmatrix} g & f\cr 1 & g \cr\end{pmatrix} \quad\hbox{ and }\quad
\hat B = \begin{pmatrix} 0 &  1 \cr  x_2\sin\theta & 0 \end{pmatrix}.
\end{equation}   
Let $M  = \begin{pmatrix} 0 & 2f x_2 \sin\theta \cr
2 & 0 \cr\end{pmatrix}$. 
Then $W_{\hat C}(\hat B) = W(M)$
is the elliptical disk with foci $\pm 2\sqrt{ fx_2\sin\theta}$
and minor axis of length $2(1-f\sin\theta |x_2|)$.
Consequently, the intersection 
of  $W_{\hat C}(\hat B)$ and the real axis always contains 
the line segment  $[f\sin\theta-1, 1-f\sin \theta]$.
\end{lemma}

\it Proof. \rm 
Since $\hat B$ is unitarily similar to 
$\hat B^t$, we have $W_{\hat C}(\hat B) = W_{\hat C}(\hat B^t)$.
By the result in \cite{Li2} (see also \cite{Li,LT}),
$W_{\hat C}(\hat B^t) = W(M)$. By Lemma \ref{ellipse}
we have the description of $W(M)$. Thus, the intersection
of $W(M)$ and the real axis has the form $[-\xi, \xi]$.
Under the assumption that $|x_2| \le 1$, 
the quantity $\xi$ will attain the minimum value $\hat \xi = 1-f\sin\theta$ 
when $x_2 = -1$ 
so that the elliptical disk $W(M)$ has foci $\pm 2i\sqrt{f\sin\theta}$
and length of minor axis $2(1-f\sin\theta)$.
Therefore, the last assertion follows.
\qed

Now, we are ready to present the following.

\medskip\noindent
{\bf Proof of Theorem \ref{main}.}
Suppose $C = \gamma I_n$ with $\gamma \in \IC$.
We may replace $C$ by $C/\gamma$ and assume that $C = I_n$.
Let $T =  E_{12}/2\in M_n$.  By Lemma \ref{cor}, 
$W_C(T) = \{0\}$.
Suppose $U$ is a unitary dilation of $T$.
By Lemma \ref{dilate}, $U$ has a compression $\hat T\in M_{n+2}$
of form (\ref{hatT}), and 
$W_C(\hat T) \subseteq W_C(U)$. 
The leading $2\times 2$ submatrix of $\hat T$ is 
$\begin{pmatrix} 0 & 1/2 \cr0 & 0 \cr \end{pmatrix}$, which
 has numerical range equal to the 
circular disk centered at the origin with radius $1/4$. So, 
for any $\xi \in \IC$ with $|\xi| \le 1/4$, there is a unitary matrix
$V = V_1 \oplus I_n$ with $V_1 \in M_2$
such that $V^*\hat TV$ has diagonal entries 
$\xi, -\xi, 0, \dots, 0, y_2/2$. 
Now, $C \oplus 0_2$ is unitarily similar to $\tilde C = [0] \oplus I_n \oplus [0]$.
Thus, $\tr(\tilde C \hat T) = -\xi\in W_C(\hat T) \subseteq W_C(U)$. 
So $W_C(U)$ contains all $\xi \in \IC$ with 
$|\xi| \le 1/4$. The conclusion of Theorem \ref{main} holds with $d = 1/4$.

\medskip
Suppose $C$ is not a scalar matrix. Let
$R = \min\{\|C-\mu I\|: \mu \in \IC\} > 0$.
By Lemma \ref{prop1}, we may 
apply a suitable unitary similarity to $C$ and  assume that 
$C$ has leading $2\times 2$ submatrix 
$\begin{pmatrix} g & f \cr R & g\cr\end{pmatrix}$
with $R = \|C-gI\| = \min\{\|C-\mu I\|: \mu \in \IC\}$. If $f = |f| e^{i\theta}$, we may
replace $C$ by $e^{-i\theta/2} D^*CD/R$ with $D = [e^{i\theta/2}]\oplus I_{n-1}$
and assume that the leading $2\times 2$ submatrix of $C$ is
$\hat C = \begin{pmatrix} g & f\cr 1 & g \cr \end{pmatrix}$. 
Since $1 = \min\{ \|C-\mu I\|: \mu \in \IC\} = \|C - gI\|$, we see that 
$f \in [0,1]$, 
the first column of $C$ has the form $(g,1,0,\dots,0)^t$ and
the second row of $C$ has the form $(1,g,0, \dots, 0)$.
We consider two cases.

\medskip\noindent
{\bf (I)} 
Suppose $1>f$. Let $T = \cos\theta E_{12} \in M_n$, where
$\theta \in [0, \pi/2)$ is sufficiently close to $\pi/2$, so that 
$$\xi_1 = 1-f\sin \theta = \cos \theta + d  \quad \hbox{ with } d > 0.$$ 
By Lemma \ref{cor}, $W_C(T)$ is a circular disk
centered at the origin with radius $\cos\theta$.
The conclusion of Theorem \ref{main} will follow once we prove the following.

\medskip
\noindent
{\bf Claim}
For any unitary dilation $U$ of $T$,
 $\xi_1 = 1-f\sin\theta = \cos \theta + d \in W_C(U)$.

\medskip
To prove the claim, let $U$ be a unitary dilation of $T$.
By Lemma \ref{dilate}, $U$ has a compression of the form $\hat T$
defined as in (\ref{hatT}). Note that
$\hat C = \begin{pmatrix} g & f \cr 1 & g\cr\end{pmatrix}$
is the leading $2\times 2$ submatrix of  $C$, 
and $\hat B = \begin{pmatrix} 0 & 1 \cr x_2 \sin \theta & 0 \cr\end{pmatrix}$ 
is the submatrix of $\hat T$ in rows and columns 
with indices 2 and $n+1$, respectively.

\smallskip
If $n = 2$, then we may assume that $C = \hat C$, and  
$W_{\hat C}(\hat B) \subseteq W_C(\hat T) \subseteq W_C(U)$.
By Lemma \ref{WCA},  $[f\sin\theta -1, 1-f\sin\theta] \subseteq W_C(U)$.
Thus,  the claim follows.

\medskip
Suppose $n \ge 3$.
Let 
$$\tilde C = 
\begin{pmatrix} g & f & c_{13} \cr 
1 & g & 0 \cr
0 & c_{32} & c_{33} \cr\end{pmatrix}, \quad  
\tilde B = \begin{pmatrix} 0 & 1 & 0 \cr x_2 \sin\theta & 0 & x_3 \cr
 0 & 0 & 0 \cr\end{pmatrix}, \quad \hbox{ and } \quad 
 B_1 = \tilde B \oplus 0_{n-3} \in M_n.$$
Then $\tilde C$ is  the leading $3\times 3$ submatrix of $C$, and  
$B_1$ 
can be obtained from  $\hat T$ by the deletion of the rows and 
columns with indices $1, n+2$, followed by a permutation similarity.
We will prove the claim by showing:
(i)  $W_{\tilde C}(\tilde B)\subseteq W_C(U)$,  
and (ii)  $\xi_1 \in W_{\tilde C}(\tilde B)$.

\medskip
To prove (i), suppose
 $\xi  =  \tr(\tilde C V_1^*\tilde B V_1)
 \in W_{\tilde C}(\tilde B)$, where  $V_1 \in M_3$
is unitary.  Let $V = V_1 \oplus I_{n-3}$. Then
$\xi = \tr(CV^*B_1V) \in W_C(\hat T) \subseteq W_C(U)$.

\medskip
To prove (ii),
let $\hat C$ and $\hat B$ be defined as in (\ref{matrixC}).
By Lemma \ref{WCA}, $\xi_1 = 1-f\sin\theta \in W_{\hat C}(\hat B)$.
So, there is  a unitary  $V_1\in M_2$ such that 
$\tr(\hat C V_1^*\hat B V_1) = \xi_1$. If $V = V_1 \oplus [\mu] \in M_3$ 
is unitary, then $\tr(\tilde C V^*\tilde B V) = \xi_1 + \bar \mu \xi_2$ with
$\xi_2 = (0,c_{32})V_1 (0, x_3)^t$.
Hence,
$$S(V_1) = \{ 
\tr (\tilde C (V_1\oplus [\mu])^*\tilde B(V_1\oplus [\mu])):
\mu \in \IC, |\mu| =  1\}$$
is a circle $S(V_1)$ with center $\xi_1$ and radius $|\xi_2|$. Now,
construct a continuous path 
of unitary matrices $V_t = e^{i(tH + (1-t)G)}\in M_2$ with $t \in [0,1]$ such that
$H, G$ are Hermitian matrices satisfying
$e^{iH} = V_1$, which is defined as above, and $V_0 = e^{iG} = 
\begin{pmatrix} 0 & 1 \cr 1 & 0 \cr\end{pmatrix}$.
Let
$$S(V_t) = \{\tr(\tilde B(V_t\oplus [\mu])^* \tilde C (V_t \oplus [\mu])):
\mu \in \IC, |\mu| =1\}\subseteq W_{\tilde C}(\tilde  B).$$
Then $S(V_t)$ is a circle with center $\xi_1(t) = 
\tr(\tilde C V_t^*\tilde B V_t)$ 
and radius $|(0, c_{32})V_t(0,x_3)^t|$.
Note that $S(V_0) = \{ x_2\sin\theta + f\}$ is a singleton.
As $t$ varies from $1$ to 0, $S(V_1)$ will change to 
the singleton $S(V_0)$ continuously. So, every point inside the circle
$S(V_1)$ will lie in some $S(V_t)$ with $t \in [0, 1]$. In particular,
$\xi_1 \in S(V_t) \subseteq W_{\tilde C}(\tilde  B)$
for some $t \in [0,1]$.
We get the desired conclusion.

\medskip\noindent
{\bf (II)}
Suppose $f = 1$, i.e., the leading $2\times 2$ submatrix of $C$ equals 
$\hat C = \begin{pmatrix} g & 1 \cr 1 & g\cr\end{pmatrix}$.
By Lemma \ref{prop1}, 
$C-gI_n$ has norm 1, which equals the $(1,2)$ entry and the $(2,1)$ entry
of the matrix. Thus, $C - gI_n$ is a direct sum of 
its leading $2\times 2$ matrix and its trailing $(n-2)\times (n-2)$ submatrix $C_1$.
Thus, $C = \hat C \oplus C_1$, and  $\|C_1-g I_{n-2}\|\le 1$.

\medskip
Let $T = E_{12}/2\in M_n$. By Lemma \ref{cor}, and the assumption on $C$,
$W_C(T)$ is a circular disk centered at 0 with radius 1/2.
We will show that there is $d > 0$ such that $1/2 + d \in W_C(U)$ 
for any unitary dilation $U$ of $T$.

Let $U$ be a unitary dilation of $T$. By Lemma \ref{dilate},
$U$ has a compression $\hat T \in M_{n+2}$ of the form (\ref{hatT}),
and $\hat T$ has a principal submatrix  
$B$ in the form (\ref{matrixB}) with $\theta = \pi/6$.
We consider two subcases.

\medskip
{\bf (II.a)} Suppose $g \notin\{1,-1\}$.
Note that $\hat T$ is permutationally similar to a matrix of the form
$$\tilde T  
= \begin{pmatrix} B &  \star \cr 0_{n-2,4} & 0_{n-2}\cr\end{pmatrix}.$$
So, for any unitary $V = V_1 \oplus I_{n-2} \in M_{n+2}$, where $V_1 \in M_4$,
$$\xi = \tr((0_2 \oplus C)V^*\tilde T V) = 
\tr((0_2 \oplus \hat C)V_1^*BV_1)$$
is an element in $W_{\hat C}(B)$.
Also, every element in $W_{\hat C}(B)$ can be put in this form.

\medskip\noindent
{\bf Claim} \it If $B$ has the form $(\ref{matrixB})$, 
then  
$W_{\hat C}(B)$  contains an interval $[0, \xi]$ with $\xi > 1/2 = \cos\theta$.
\rm

\medskip
Suppose the claim holds.
Since the set  $\{(x_1, x_2, y_1, y_2)^t \in \IC^4: 
|x_1|^2 + |x_2|^2 \le 1, |y_1|^2 + |y_2|^2 \le 1\}$  is compact,
there is $d > 0$ such that  $W_{\hat C}(B)$ will contain a number larger 
than $\cos\theta+d = 1/2 + d$ for any matrix $B$ in the form (\ref{matrixB}).
Consequently, $1/2+d \in W_C(U)$ for any unitary dilation $U$ of $T$.

\medskip
To prove the claim, 
let  $V = \frac{1}{\sqrt 2} \begin{pmatrix} 1 & 1 \cr 1 & -1 \cr\end{pmatrix} 
\oplus I_2$,
and  $(c,s)  = (1/2,\sqrt 3/2)$. Then 
$$\tilde B = V^*BV = \begin{pmatrix}
c/2   & -c/2 & 1/\sqrt{2} & -s/\sqrt{2}\cr
c/2   & -c/2 & -1/\sqrt{2} & -s/\sqrt{2}\cr
x_1/\sqrt{2} + x_2s/\sqrt{2} & x_1/\sqrt{2} - x_2s/\sqrt{2} & 0 & x_2c\cr
y_1/\sqrt{2} + y_2s/\sqrt{2}& y_1/\sqrt{2} - y_2s/\sqrt{2} & 0 & y_2c \cr
\end{pmatrix}.$$
For $j = 1,2$, let $B_j \in M_3$ be obtained from $\tilde B$ by deleting the
$j$th row and $j$th column so that
$$B_1 =  \begin{pmatrix}
 -c/2 & -1/\sqrt{2} & -s/\sqrt{2}\cr
 x_1/\sqrt{2} - x_2s/\sqrt{2} & 0 & x_2c\cr
 y_1/\sqrt{2} - y_2s/\sqrt{2} & 0 & y_2c \cr
\end{pmatrix} \ \hbox{ and } \ 
 B_2 = \begin{pmatrix}
c/2    & 1/\sqrt{2} & -s/\sqrt{2}\cr
x_1/\sqrt{2} + x_2s/\sqrt{2} & 0 & x_2c\cr
y_1/\sqrt{2} + y_2s/\sqrt{2} &  0 & y_2c \cr
\end{pmatrix}.$$
We consider two cases.

\medskip
{\bf Case 1} Suppose $-c/2$ is a boundary point of $W(B_1)$ and 
$c/2$ is a boundary point of $W(B_2)$.
By Lemma \ref{lemma}, 
$$1 = |x_1 + x_2 s| = |x_1 - x_2 s| \quad \hbox{ and } \quad 
s = |y_1 + y_2 s| = |y_1 - y_2 s|.$$
Thus, $0 =  x_1 \bar x_2 + \bar x_1 x_2, 1 = |x_1|^2 + s^2|x_2|^2$.
Since $|x_1|^2 + |x_2|^2 \le 1, |y_1|^2 + |y_2^2| \le 1$, we see that
$x_2 = 0$ and $|x_1| = 1$.
Since $|x_1|^2 + |y_1|^2 \le 1$, it follows that $y_1 = 0$, and
 $|y_2s| = s$, i.e., $|y_2| = 1$. 
But then the matrix $\hat B$ in (\ref{matrixC})
will be of the form 
$\begin{pmatrix} 0 & 1 \cr 0 & 0 \cr\end{pmatrix}$
so that $W_{\hat C}(\hat B) = W_{\hat C}(\hat B^t) = W(M)$,
where $M = \begin{pmatrix} 0 & 0 \cr 2 & 0 \cr\end{pmatrix}$ by
Lemma \ref{WCA}. Thus,  $W_{\hat C}(\hat B)$
 is the unit disk containing the interval $[0,1]$.

\medskip
{\bf Case 2} Suppose
$-c/2$ is an interior point of $W(B_1)$ or 
$c/2$ is an interior point of $W(B_2)$.
Here recall that we assume that $g \notin \{1, -1\}$.
If $-c/2$ is an interior point of $B_1$, then there is $\delta > 0$ such that
$-c/2 + \varepsilon/(g-1) \in W(B_1)$ for any $|\varepsilon|< \delta$, and 
$\tilde B$ is unitarily similar to a matrix with its leading 
$2\times 2$ submatrix 
$$B_3 = \begin{pmatrix} c/2 & \star \cr 
\star & -c/2+\varepsilon/(g-1)\cr\end{pmatrix}.$$
The matrix $\hat C = \begin{pmatrix} g& 1 \cr 1 & g \cr\end{pmatrix}$
is unitarily similar to 
$\tilde C = \diag(g+1, g-1)$.
Thus, $W_{\tilde C}(B_3)$ contains real numbers of the form
$(g+1)c/2 + (g-1)(-c/2 + \varepsilon/(g-1)) = c + \varepsilon$.
Hence,  $W_{\tilde C}(B_3)$ contains the interval $[c, c+\varepsilon]$.

\medskip
Similarly, if $c/2$ is an interior point of $B_2$, then there is 
$\delta > 0$ such that
$c/2 + \varepsilon/(g+1) \in W(B_2)$ whenever $|\varepsilon| < \delta$, and 
$\tilde B$ is unitarily similar to a matrix with its leading 
$2\times 2$ submatrix 
$$B_4 = \begin{pmatrix} c/2+\varepsilon/(g+1) & \star \cr 
\star & -c/2\cr\end{pmatrix}.$$
Then, $W_{\tilde C}(B_4)$ contains  real numbers of the form
$(g+1)(c/2 + \varepsilon/(g+1)) - (g-1)c/2 = c + \varepsilon$. 
Hence,  $W_{\tilde C}(B_4)$ contains the interval $[c, c+\varepsilon]$.

\medskip
Combining Case 1 and Case 2, we establish the  claim. The  theorem follows.

\medskip
{\bf (II.b)} Suppose $g \in \{1, -1\}$. Then $\hat C$ is  unitarily similar to 
$\diag(2g,0)$.
We may assume that $g = 1$. Otherwise, replace $C$ by $-C$.
Thus, $W_{\hat C}(B) = 2W(B)$.
Since $C$ is not a rank one normal matrix, $C_1 \ne 0$. 
Thus, $C_1$ is unitarily similar to a matrix with a nonzero $(1,1)$ entry
$\mu$.

{\bf Case 1.}
Suppose in the matrix $B$,
$|x_2| \le \sqrt 3/6$. By Lemma \ref{ellipse}, the submatrix 
$B_0 = \begin{pmatrix} 0 & 1 \cr x_2 \sqrt{3}/2 & 0 \cr
\end{pmatrix}$ of $B$ has numerical range equal to an elliptical 
with the length of minor axis  $1 - |x_2|\sqrt 3/6 \ge 1 - 1/4 = 3/4$. Thus,
there is a unitary matrix $V = [1] \oplus V_1 \oplus [1]$ with $V_1 \in M_2$
such that $V^*BV$ has diagonal entries
$ 0, \xi,  -\xi, y_2/2$ for any $|\xi| \le \sqrt 3/4$, and 
$\hat T$ is unitarily similar to a matrix of the form
$\hat T_\xi = \begin{pmatrix} 
V^*BV & \star  \cr 0_{n-2,4} & 0_{n-2} \cr \end{pmatrix}$.
Since  $C\oplus 0_2$ is unitarily similar to 
$\tilde C = \diag(0,2,0,0) \oplus C_1$, 
$\tr(\tilde C \hat T_\xi) = \xi.$ So, $W_C(U)$ always contains
$\xi$ with $|\xi| \le \sqrt 3/4$.

{\bf Case 2.} If $|x_2| \ge \sqrt 3/6$, then  the submatrix 
$\hat B_0 = \begin{pmatrix} 0 & x_2/2 \cr 0 & y_2 /2\cr
\end{pmatrix}$ of $B$ has numerical range equal to an elliptical 
with foci $0, y_2/2$ and length of minor axis 
$|x_2|/2 \ge \sqrt 3/12$. Thus, the focus 0 is an interior point of 
$W(\hat B_0)$, and there is $\delta > 0$ such that
$\xi \in W(\hat B_0)$ whenever $|\xi| < \delta$.
As a result, for any $\xi_1, \xi_2$ with  $|\xi_1| \le 1/4$ and $|\xi_2| \le \delta$,
there is a unitary matrix $V = V_1 \oplus V_2$ with $V_1, V_2\in M_2$
such that $V^*BV$ has diagonal entries
$\xi_1,-\xi_1,  y_2/2 -\xi_2, \xi_2$,
and $\hat T$ is unitarily similar to a matrix of the form
$\hat T_\xi = \begin{pmatrix} 
V^*BV & \star \cr 0_{n-2,4} & 0_{n-2} \cr
\end{pmatrix}$.

If $n = 3$, then $C_1 = [\mu]$, and 
$C\oplus 0_2$ is unitarily similar to 
$\tilde C = \diag(2,0,0,\mu, 0)$. We have 
$\tr(\tilde C \hat T_\xi) = \xi_1 + \mu \xi_2$.
Hence, 
$W_C(U)$ always contains $\xi$ with $\xi \in (0, 1/4 + |\mu|\delta]$.

\medskip
If $n \ge 4$, we may assume that  $C\oplus 0_2$ is unitarily similar to
$\tilde C = \diag(2,0,0,0) \oplus C_1$.
Suppose $\tilde C$ has columns $\bc_1, \dots, \bc_{n+2}$.
Then $\bc_5 = (0,0,0,0, \mu, \eta_1, \dots, \eta_{n-3})^t$ 
with $\eta_1, \dots, \eta_{n-3}\in \IC$.
We may replace $\tilde C$ by 
$(I_5 \oplus P_1)^* \tilde C (I_5 \oplus P_1)$, 
where $P_1 \in M_{n-3}$ is a unitary matrix such that 
$(I_5\oplus P_1)^* \bc_5 = (0,0,0,0,\mu, \eta, 0, \dots, 0)^t$
with $\eta = \sqrt{\sum_{j=1}^{n-3} |\eta_j|^2}$.
Note that this change does not affect the leading $5\times 5$ submatrix 
of the original $\tilde C$.

\medskip
Suppose $\hat T_\xi$ have rows $\bd_1, \dots, \bd_{n+2}$. 
Let $\bd_4 = (d_{41}, \dots, d_{4,n+2})$. Then  $d_{44} =\xi_2$. 
We may replace  $\hat T_\xi$ by
$(I_4 \oplus P_2)^*\hat T_\xi(I_4 \oplus P_2)$,
where $P_2 \in M_{n-2}$ is a unitary matrix such that
$\bd_4 (I_4 \oplus P_2) = (d_{41},d_{42},d_{43}, \xi_2, z, 0, \dots, 0)$
with $z = \sqrt{\sum_{j=5}^{n+2}|d_{4j}|^2}$.
Note that this change does not affect the leading $4\times 4$ submatrix of 
the original $\hat T$.

\medskip
Now, let $\tilde T_\xi$ be obtained from $\hat T_\xi$ by switching rows 4 and 5,
and also switching columns 4 and 5.
Suppose 
$\tilde T_\xi$ have rows $\bg_1, \dots, \bg_{n+2}$. 
Then $\bg_1, \bg_2, \bg_3, \bg_5$ are the only  nonzero rows, which are obtained from 
$\bd_1, \bd_2, \bd_3, \bd_4$ defined in the preceding paragraph
by switching the 4$th$ and 5$th$ entries in the
vectors.
Since $\bg_1 \bc_1 = \xi_1$, $\bg_5 \bc_5 = 
(d_{41}, d_{42}, d_{43}, z, \xi_2, 0, \dots, 0) (0,0,0,0, \mu, \eta, 0, \dots, 0)^t
= \mu\xi_2$, $\bc_2 = \bc_3 = \bc_4 = 0_{n+2}$,
it follows that
$$\tr(\tilde C \tilde T_\xi) = 
\tr(\tilde T_\xi\tilde C) = \sum_{j=1}^{n+2} \bg_j \bc_j = 
\bg_1 \bc_1 + \bg_5 \bc_5 = \xi_1 + \mu \xi_2.$$
So, $W_C(U)$ always contains $\xi$ with $\xi \in (0, 1/4 + |\mu|\delta]$.

\medskip
Let $1/2 + d = \min\{\sqrt 3/4, 1/4 + |\mu|\delta\}$. Then
$W_C(U)$ always contains $1/2+d$ by the analysis in Case 1 and Case 2.
The  theorem follows.
\qed


\medskip
\noindent
{\bf Proof of Theorem \ref{key}}
The implication (a) $\Rightarrow$ (b) follows from the result
in \cite{CL}. The implication (b) $\Rightarrow$ (c) is clear.
By Theorem \ref{main},  we have the implication
(c) $\Rightarrow$ (a). \qed

\medskip\noindent
{\bf Acknowledgment} 

The author would like to thank the referee for 
her/his careful reading of the manuscript, and helpful comments.

\medskip\noindent
{\bf Declaration.} 

There is no conflict of interest connected to this article.

\end{document}